\documentclass[MSNbibl,number,citesort,seceqn,dvips]{arxbj}
\usepackage{upgreek}
\usepackage{graphicx}
\aid{0}
\volume{19}
\issue{5B}
\pubyear{2013}
\firstpage{2689}
\lastpage{2714}
\doi{10.3150/12-BEJ471} 

\makeatletter
\newcommand{\RMo}{\mathrm{o}}
\newcommand{\RMO}{\mathrm{O}}

\newcommand{\mrmd}{\,\mathrm{d}}
\newcommand{\mrmdd}{\mathrm{d}}

\newcommand{\bigtimes}{\mathop{\mbox{\fontsize{17}{17}\selectfont{$\!\times$}}}}

\newtheorem{prop}{Proposition}

\newtheorem{property}{Property}

\newremark{example}{Example}

\newcommand{\Prob}{\mbox{\textsf{P}}}
\newcommand{\Ex}{\mathbf{X}_E}
\newcommand{\Pa}{\mathbf{X}_P}
\newcommand{\Lsv}{\mathcal{L}}

\makeatother

\begin{document}
\begin{frontmatter}

\title{A new representation for multivariate tail~probabilities}
\runtitle{Multivariate tail probabilities}

\begin{aug}
\author[1]{\fnms{J.L.} \snm{Wadsworth}\corref{}\thanksref{1}\ead[label=e1]{jenny.wadsworth@epfl.ch}} \and
\author[2]{\fnms{J.A.} \snm{Tawn}\thanksref{2}\ead[label=e2]{j.tawn@lancaster.ac.uk}}
\runauthor{J.L. Wadsworth and J.A. Tawn} 
\address[1]{Ecole Polytechnique F\'{e}d\'{e}rale de Lausanne, EPFL SB
MATHAA STAT,
Station 8, 1015 Lausanne, Switzerland. \printead{e1}}
\address[2]{Department of Mathematics and Statistics, Lancaster
University, LA1 4YF,
United Kingdom.\\ \printead{e2}}
\end{aug}

\received{\smonth{10} \syear{2011}}
\revised{\smonth{7} \syear{2012}}

%
\begin{abstract}
Existing theory for multivariate extreme values focuses upon
characterizations of the distributional tails when all components of a
random vector, standardized to identical margins, grow at the same
rate. In this paper, we consider the effect of allowing the components
to grow at different rates, and characterize the link between these
marginal growth rates and the multivariate tail probability decay rate.
Our approach leads to a whole class of univariate regular variation
conditions, in place of the single but multivariate regular variation
conditions that underpin the current theories. These conditions are
indexed by a homogeneous function and an angular dependence function,
which, for asymptotically independent random vectors, mirror the role
played by the exponent measure and Pickands' dependence function in
classical multivariate extremes. We additionally offer an inferential
approach to joint survivor probability estimation. The key feature of
our methodology is that extreme set probabilities can be estimated by
extrapolating upon rays emanating from the origin when the margins of
the variables are exponential. This offers an appreciable improvement
over existing techniques where extrapolation in exponential margins is
upon lines parallel to the diagonal.
\end{abstract}

%

\end{frontmatter}

\section{Introduction}
\label{secIntroduction}

The concept of regular variation underpins multivariate extreme value
theory. However, in recent years the understanding of multivariate
extremes has broadened and a number of different regular variation
conditions have emerged as useful tools for analyzing extremal
dependence. Classical multivariate extreme value theory (e.g.,
de Haan and de Ronde~\cite{deHaandeRonde98}) provides a helpful
characterization of the
limit distribution of suitably normalized random vectors when all
components of the vector are asymptotically dependent. We focus
initially on the bivariate case. For any bivariate random vector
$(X,Y)$ with marginal distribution functions $F_{X},F_{Y}$, the notion
of asymptotic dependence is defined through the limit
%
\begin{equation}\label{chi}
\chi:=\lim_{u\rightarrow1} \Prob\bigl\{F_{X}(X)>u|F_{Y}(Y)>u
\bigr\}.
\end{equation}
The cases $\chi>0$ and $\chi=0$ define, respectively, \textit{asymptotic
dependence} and \textit{asymptotic independence} of the vector. It is
convenient to standardize the margins, which is theoretically achieved
by the probability integral transform; we assume continuous marginal
distribution functions. Let $(X_P,Y_P)$ be a random vector with
standard Pareto marginal distributions, that is, $\Prob(X_P>x)= \Prob
(Y_P>x)= 1/x, x\geq1$. The classical bivariate extreme value paradigm
is characterized by the following convergence of measures. We say that
$(X_P,Y_P)$ is in the \textit{domain of attraction} of a bivariate
extreme value distribution if there exists a measure $\mu$ such that
for every relatively compact Borel set $B \subset[0,\infty
]^2\setminus
\{\mathbf{0}\}$,
%
\begin{equation}\label{MSconv}
\lim_{n\rightarrow\infty} n\Prob\bigl\{(X_P/n,Y_P/n) \in B
\bigr\} = \mu(B),
\end{equation}
where $\mu$ is non-degenerate, and $\mu(\partial B)=0$. Equivalently
one can write $n\Prob\{(X_P/n,\break Y_P/n) \in\cdot\}\stackrel
{v}{\rightarrow
} \mu(\cdot)$, the symbol $\stackrel{v}{\rightarrow}$ denoting vague
convergence of measures (see, e.g., Resnick~\cite{Resnick06}). The limit
measure $\mu$ is homogeneous of order $-1$ and has mass on the subcone
$(0,\infty]^2 \subset[0,\infty]^2 \setminus\{\mathbf{0}\}$ only if the
vector exhibits asymptotic dependence. Under asymptotic independence,
the mass of $\mu$ concentrates upon the axes of $[0,\infty]^2
\setminus
\{\mathbf{0}\}$; equivalently the limiting joint distribution of normalized
componentwise maxima of $(X_P,Y_P)$ corresponds to independence.

From a practical perspective, the major shortcoming of
representation (\ref{MSconv}) is that variables can exhibit marked
dependence at all observable levels, yet be asymptotically independent.
Consequently, in this case, the use of limit models of the form (\ref
{MSconv}) for the dependence is inappropriate for extrapolation from
the observed data to more extreme levels.

Ledford and Tawn \cite{LedfordTawn96,LedfordTawn97} addressed this
problem by
constructing a theory for the rate of decay of the dependence under
asymptotic independence. They introduced the concept of the \textit
{coefficient of tail dependence}, $\eta\in(0,1]$, which characterizes
the order of decay of $\Prob\{(X_P/n,Y_P/n) \in B\}$ for any Borel set
$B \subset(0,\infty]^2$. The theory has since been elaborated on in the
works of Resnick and co-authors (Resnick \cite{Resnick02},
Maulik and Resnick \cite{MaulikResnick04}, Heffernan and Resnick \cite
{HeffernanResnick05}), under the label
\textit{hidden regular variation}. Under asymptotic independence,
$\Prob\{ (X_P/n,Y_P/n) \in B\} = \RMo(n^{-1})$ for $B
\subset(0,\infty]^2$. This rate can be adjusted to $\RMO(n^{-1})$ by
altering the normalizations. The random vector $(X_P,Y_P)$ possesses
hidden regular variation if both limit (\ref{MSconv}) holds, and in
addition on $(0,\infty]^2$,
%
\begin{equation}\label{HRVconv}
n\Prob\bigl\{\bigl(X_P/c(n),Y_P/c(n)\bigr) \in\cdot
\bigr\} \stackrel{v} {\rightarrow} \nu(\cdot)
\end{equation}
with $\nu$ non-degenerate, and $c(n)=\RMo(n)$, $n\to\infty$. In
representation (\ref{HRVconv}), the sequence $c(n)$, is a regularly
varying function of $n$ with index $\eta$, whilst the limit measure
$\nu
$ is homogeneous of order $-1/\eta$ (Resnick \cite
{Resnick02,Resnick06}).

Ledford and Tawn \cite{LedfordTawn97} phrased their original condition
in terms of the
particular set $B=(x,\infty]\times(y,\infty]$. Their assumption, in
Pareto margins, was that as $n\rightarrow\infty$,
%
\begin{equation}\label{LT}
\Prob(X_P>nx,Y_P>ny) \sim n^{-1/\eta}x^{-c_1}y^{-c_2}
\Lsv(nx,ny),
\end{equation}
where, $c_1,c_2\geq0$, $c_1+c_2=1/\eta$ and $\Lsv$ is a bivariate
slowly varying function. That is, $\Lsv$ satisfies $\lim_{n\rightarrow
\infty}\Lsv(nx,ny)/\Lsv(n,n) = g(x,y)$, $x,y>0$, with $g$ homogeneous
of order 0. Consequently, the value of $g$ depends only upon the ray
$w:=x/(x+y)$ and so $g^*(w):= g(w, 1-w)$ is defined as the \textit{ray
dependence function}. For identifiability of all components in
representation~(\ref{LT}), $g^*$ also satisfies a quasi-symmetry
condition; see Ledford and Tawn \cite{LedfordTawn97}. The formulation
is discussed in
greater detail in Section \ref{secConnectionsToExistingTheory}.

In both representations (\ref{MSconv}) and (\ref{HRVconv}), the key
idea behind the limit theory is that the rate of tail decay of the
bivariate distribution of $(X_P/n,Y_P/n)$ is uniform across all rays
emanating from the origin. The rates are $\RMO(n^{-1})$ under asymptotic
dependence and $\mathrm{O}\{n^{-1/\eta} \Lsv(n,n)\}$ under asymptotic
independence. Let $(X_E,Y_E)=(\log X_P, \log Y_P)$ be a random vector
with standard exponential margins. Statistical methodology for
inference upon sets which are extreme in all variables is based on the
equivalent representations
%
\begin{eqnarray}\label{LText}
\Prob\bigl\{(X_P,Y_P) \in tnB\bigr\} &\sim&
t^{-1/\eta}\Prob\bigl\{(X_P,Y_P) \in nB\bigr\},
\nonumber\\[-8pt]\\[-8pt]
\Prob\bigl\{(X_E,Y_E) \in v + \log n + A\bigr\} &\sim&
\exp(-v/\eta) \Prob\bigl\{ (X_E,Y_E) \in\log n+ A\bigr
\}\nonumber
\end{eqnarray}
as $n \rightarrow\infty$, for Borel sets $B \subset(0,\infty]^2$, $A
\subset(-\infty,\infty]^2$, $t\geq1, v\geq0$, with addition applied
componentwise. Practically, $\eta$ is estimated using one of several
methods (see, e.g., Beirlant \textit{et al.}~\cite{Beirlant04},
Chapter 9, for a review of
available methodology), whilst\vspace*{1pt} the sets $nB$, $\log n + A$ become some
extreme sets $B' \subseteq[1,\infty]^2$ and $A'\subseteq[0,\infty]^2$.
The probabilities $\Prob\{(X_P,Y_P) \in B'\}$, $\Prob\{(X_E,Y_E) \in
A'\}$ are then estimated empirically. There is a trade-off between the
sets $B', A'$ being sufficiently extreme for the asymptotic
representation (\ref{LText}) to be reliable, and containing
sufficiently many points for the empirical probability estimate to be
reliable. Under this theory, one can extrapolate upon rays $y=\{
(1-w)/w\}x$, $w\in(0,1)$, when the margins are Pareto, or lines of the
form $y=x+d, d=\log\{(1-w)/w\}$, when the margins are exponential.

When asymptotic independence is present, we may not be interested in
events where all variables are simultaneously extreme.
Heffernan and Tawn \cite{HeffernanTawn04} presented a new limit theory
to address this problem.
For a random vector $(X_G,Y_G)$ with standard Gumbel marginal
distributions, they considered the limiting distribution of the
variable $X_G$, appropriately affinely normalized, conditional upon the
concomitant variable $Y_G$ being extreme. This again can be phrased as
another regular variation condition upon the cone $[0,\infty]\times
(0,\infty]$; see Heffernan and Resnick~\cite{HeffernanResnick07}.
Specifically, the
assumption in Heffernan and Resnick \cite{HeffernanResnick07} can be
expressed in
exponential margins as there exist functions $s(n)>0, r(n)$ such that
\[
\Prob\bigl(\bigl\{X_E - r(n)\bigr\}/s(n) \leq x|Y_E >
\log n \bigr) \rightarrow H(x)
\]
as $n\rightarrow\infty$ for some non-degenerate distribution function
$H$. Under the methodology of Heffernan and Tawn \cite
{HeffernanTawn04}, the
normalization functions $r,s$ are estimated parametrically, whilst the
distribution function $H$ is estimated empirically.

\begin{figure}

\includegraphics{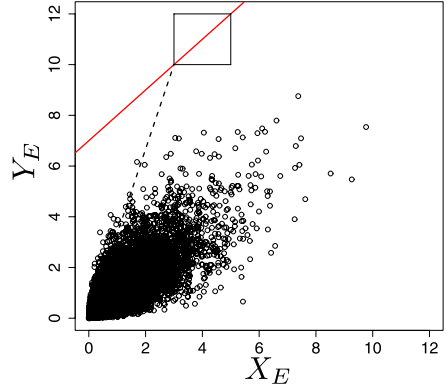}

\caption{Bivariate normal data with dependence parameter $\rho
=0.8$ on exponential margins. The line parallel to the diagonal
represents a trajectory for extrapolation under the existing
Ledford and Tawn \cite{LedfordTawn97} methodology; the dashed line
represents a trajectory
for extrapolation based upon the methodology we outline in Section
\protect\ref{secStatisticalMethodology}.}
\label{figNorm08Exp}
\end{figure}

In the theory of hidden regular variation, the normalization function
$c(n)$, defined by limit (\ref{HRVconv}), is determined by the
extremal dependence structure of $(X_P,Y_P)$ upon fixing $\Prob\{
(X_P/c(n),Y_P/c(n))\in B\} \asymp n^{-1}$. In this paper, we explore
the effects of applying different normalizations to the components and
characterizing the decay rates of the associated multivariate
distribution. In particular, we explore the behaviour of
%
\begin{equation}\label{BG}
\Prob\bigl(X_P>n^\beta,Y_P>n^\gamma
\bigr) \equiv\Prob(X_E>\beta\log n,Y_E>\gamma\log n)
\end{equation}
for $(\beta,\gamma)\in\mathbb{R}^2_+\setminus\{\mathbf{0}\}$ as
$n\rightarrow\infty$. This permits a generalization of the theory of
Ledford and Tawn \cite{LedfordTawn97}, which is a special case with
the restriction
$\beta=\gamma$. This generalization is desirable as it promotes
methodology which overcomes some of the shortcomings of
Ledford and Tawn~\cite{LedfordTawn97}, of which the principal drawback
is the limited
trajectories in which one can extrapolate. To highlight this, consider
Figure \ref{figNorm08Exp}, which shows realizations from a standard
bivariate normal with dependence parameter $\rho=0.8$ on exponential
margins. The solid line in the picture indicates a trajectory in which,
for certain sets such as the one indicated, one would not be able to
use the methodology of Ledford and Tawn \cite{LedfordTawn97}, with
empirical estimation
of $\Prob\{(X_P,Y_P) \in A'\}$, to yield any non-zero estimation of a
probability. A similar problem can be incurred in the methodology of
Heffernan and Tawn \cite{HeffernanTawn04}, as identified by
Peng and Qi \cite{PengQi04}, see
Section \ref{secStatisticalMethodology}. Ramos and Ledford \cite
{RamosLedford09} offer
a parametric approach following on from Ledford and Tawn \cite
{LedfordTawn97} which can
avoid this issue. However, the use of parametric models brings concerns
of mis-specification, and furthermore is only valid as an approximation
on a sub-region $(u,\infty)^2$ for large $u$. The theory and subsequent
methodology that we develop permits extrapolation upon rays in
exponential margins, illustrated by the dashed line in Figure \ref
{figNorm08Exp}. As such, this substantially broadens the range of
sets for which non-zero estimation is possible, whilst avoiding
parametric assumptions. Our new representation, although more closely
tied with the Ledford and Tawn \cite{LedfordTawn97} paradigm, also
provides some modest
links to the alternative limit theory of Heffernan and Tawn \cite
{HeffernanTawn04} and
Heffernan and Resnick \cite{HeffernanResnick07}, and permits
extrapolation into regions
where not all components are simultaneously extreme.

The motivating philosophy behind our approach is a little different to
the standard form. The usual limit theory for multivariate extremes
takes a single normalization, and characterizes the `interesting'
(i.e., non-degenerate) limits that arise. By contrast here, we
examine a whole class of different normalizations, but initially place
our focus upon the asymptotic decay rate associated with each
normalization. This characterization of the links between
normalizations and decay rates enriches the mathematical structure
behind multivariate extreme value theory and promotes novel statistical
methodologies.

In Section \ref{secCharacterizationOfDecayRates}, we develop the new
representation and its justification. We focus here upon the link
between the normalizing functions and the tail probability decay rates
of the distributions. We also consider extensions to higher dimensions
and associated issues. Section \ref{secConditionalMultivariateLimits}
explores extensions to bivariate limits. In Section \ref
{secConnectionsToExistingTheory}, we detail more fully the links to
the existing theories for asymptotically independent bivariate
extremes. We provide some specific examples covering a variety of
dependence structures, detailing their characterizations in each of the
different extreme value models. Some possible statistical methodology
is suggested and illustrated in Section \ref
{secStatisticalMethodology}; a~concluding discussion is provided in
Section \ref{secDiscussion}.


\section{Characterization of decay rates}
\label{secCharacterizationOfDecayRates}

\subsection{Regular variation assumption}
\label{secBasicAssumption}

We focus on the characterization of the leading order power behaviour
of expression (\ref{BG}) as \mbox{$n\rightarrow\infty$}. Following similar
assumptions to Ledford and Tawn \cite{LedfordTawn97}, we suppose that
for all $(\beta,\gamma)\in\mathbb{R}^2_+\setminus\{\mathbf{0}\}$,
%
\begin{equation}\label{A1}
\Prob\bigl(X_P>n^\beta,Y_P>n^\gamma
\bigr) =\Prob(X_E>\beta\log n,Y_E>\gamma\log n) = L(n;
\beta,\gamma)n^{-\kappa(\beta,\gamma)},
\end{equation}
where $L$ is a univariate slowly varying function in $n$, $n\to\infty$,
for all $(\beta,\gamma)\in\mathbb{R}^2_+\setminus\{\mathbf{0}\}$, and the
function $\kappa(\beta,\gamma) > 0$ maps the different marginal growth
rates to the joint tail decay rate. Equation (\ref{A1}) implies that
$\Prob(X_P>n^\beta,Y_P>n^\gamma)$ is a regularly varying function at
infinity of index $-\kappa(\beta,\gamma)$.


\subsection{Consequences of the assumption}
\label{secBasicResults}

\subsubsection*{The index of regular variation}
\label{secTheIndexOfRegularVariation}

The function $\kappa$ is the key quantity in determining the behaviour
of expression (\ref{A1}) as $n\rightarrow\infty$. We now explore some
of its properties. By the regular variation assumption (\ref{A1}), it
follows (e.g., Resnick \cite{Resnick06}, Proposition 2.6(i)) that
%
\begin{equation}\label{kappaRV}
\kappa(\beta,\gamma) = -\lim_{n\rightarrow\infty}\log\Prob
\bigl(X_P>n^\beta,Y_P>n^\gamma\bigr)/\log n.
\end{equation}

\begin{property}
\label{property1}
$\kappa$ is homogeneous of order 1.
\end{property}
\begin{pf}
Take $h>0$. Then
\begin{eqnarray*}
\kappa(h\beta,h\gamma) &=& -\lim_{n\rightarrow\infty}\log\Prob\bigl
(X_P>n^{h\beta},Y_P>n^{h\gamma}
\bigr)/\log n
\\
&=& -h\lim_{n\rightarrow\infty}\log\Prob\bigl(X_P>n^{h\beta
},Y_P>n^{h\gamma
}
\bigr)/(h\log n)
\\
&=& h\kappa(\beta,\gamma).
\end{eqnarray*}
\upqed\end{pf}
%
\begin{property}
\label{property2}
$\kappa$ is a non-decreasing function of each argument.
\end{property}
\begin{pf}
Let $t>1$, thus $\Prob(X_P>n^{t\beta},Y_P>n^\gamma) \leq\Prob
(X_P>n^\beta,Y_P>n^\gamma)$. Following this inequality through
limit (\ref{kappaRV}) yields $\kappa(t\beta,\gamma)\geq\kappa
(\beta,\gamma)$.
\end{pf}
%
\begin{property}
\label{property3}
$\kappa(\beta,0) = \beta; \kappa(0,\gamma) = \gamma$.
\end{property}
This holds by the assumptions on the margins, that is, $\Prob
(X_P>n^{\beta}, Y_P>1) = \Prob(X_P>n^{\beta}) = n^{-\beta}$.
%
\begin{property}
\label{property4}
For any dependence structure, $\max\{\beta,\gamma\} \leq\kappa
(\beta,\gamma)$. Under positive quadrant dependence of $(X_P, Y_P)$,
$\kappa
(\beta,\gamma) \leq\beta+\gamma$.
\end{property}
\begin{pf}
Complete dependence provides the bound $\Prob(X_P>n^\beta,
Y_P>n^\gamma
) \leq\Prob(X_P> \max\{n^\beta,\allowbreak n^\gamma\})$, whilst positive quadrant
dependence implies $\Prob(X_P>n^\beta)\Prob(Y_P>n^\gamma) \leq
\Prob
(X_P>n^\beta, Y_P>n^\gamma)$. Following these inequalities through
limit (\ref{kappaRV}) yields the result.
\end{pf}

Note that without positive quadrant dependence, $\kappa(\beta,\gamma)$
can become arbitrarily large; the bivariate normal distribution with
negative dependence (Example \ref{secBivariateNormalDistribution} of
Section \ref{secExamples}) provides an example of this.

The function $\kappa$ provides information about the level of
dependence between variables at sub-asymptotic levels. The homogeneity
property of $\kappa$ suggests it is instructive to consider a
decomposition into a radial and angular component. Define pseudo-angles
$\omega:=\beta/(\beta+\gamma)$. Then $\kappa(\beta,\gamma) =
(\beta
+\gamma)\kappa(\omega,1-\omega)$, which motivates defining an angular
dependence function $\lambda(\omega):=\kappa(\omega,1-\omega),
\omega\in
[0,1]$. Setting $\omega=1/2$ returns the coefficient of tail dependence
$\eta$ as defined in Ledford and Tawn \cite{LedfordTawn96}, that is,
$\kappa
(1,1)=2\lambda(1/2)=1/\eta$. Interestingly, the link between $1/\eta$
and the newly-defined function $\lambda$ is analogous to the link
between the extremal coefficient $\theta$ and Pickands' dependence
function $A$ (Pickands \cite{Pickands81}) of classical multivariate extremes:
$\theta=2A(1/2)$, see, for example, Beirlant \textit{et al.} \cite
{Beirlant04}.

In terms of $\lambda$, the inequalities of Property \ref{property4}
read: $\max\{\omega,1-\omega\} \leq\lambda(\omega) \leq1$. This
further highlights the parallels between $\lambda$ and Pickands'
dependence function $A$, since this is one of the two characterizing
features of $A$. The other feature of $A$ is convexity; convexity of
$\kappa$ and $\lambda$ entails an additional dependence condition.
%
\begin{prop}
\label{prop3}
For a random vector satisfying assumption (\ref{A1}), the homogeneous
function $\kappa$ and angular dependence function $\lambda$ are convex
if there exists $N \in\mathbb{N}$ such that for all $n>N$, and any
$\delta_1, \delta_2, \varepsilon_1, \varepsilon_2 \geq0$, $\delta_1+\varepsilon
_1> 0$ and $\delta_2+\varepsilon_2 > 0$,
%
\begin{eqnarray}\label{ineq}
&&
\Prob\bigl\{X_{E}> (\delta_1 + \delta_2)
\log n, Y_{E}> (\varepsilon_1 + \varepsilon_2) \log
n \bigr\} \nonumber\\[-8pt]\\[-8pt]
&&\quad\geq\Prob(X_{E}> \delta_1 \log n,
Y_{E}> \varepsilon_1 \log n ) \Prob(X_{E}>
\delta_2 \log n, Y_{E}> \varepsilon_2 \log n ).\nonumber
\end{eqnarray}
Conversely if $\kappa, \lambda$ are strictly convex, then for each such
$\delta_1, \delta_2, \varepsilon_1, \varepsilon_2$, there exists
$N(\delta_1,\delta_2,\allowbreak\varepsilon_1,\varepsilon_2) \in\mathbb{N}$ such that
for all
$n>N(\delta_1,\delta_2,\varepsilon_1,\varepsilon_2)$, inequality (\ref
{ineq}) holds.
\end{prop}
\begin{pf}
Suppose inequality (\ref{ineq}) holds. By taking log of each side
of (\ref{ineq}), and applying (\ref{kappaRV}), we deduce $\kappa
(\delta_1 + \delta_2, \varepsilon_1 + \varepsilon_2) \leq\kappa(\delta_1,
\varepsilon_1) + \kappa(\delta_2, \varepsilon_2)$, that is, $\kappa$ is subadditive.
Homogeneous order 1 functions defined on convex cones are convex if and
only if they are subadditive (e.g., Niculescu and Persson \cite
{NiculescuPersson06}, Lemma
3.6.1). In terms of $\lambda$, we have
\[
(\delta_1 + \delta_2 + \varepsilon_1 +
\varepsilon_2) \lambda\biggl(\frac
{\delta_1 + \delta_2}{\delta_1 + \delta_2 + \varepsilon_1 + \varepsilon_2}
\biggr) \leq(
\delta_1 + \varepsilon_1) \lambda\biggl(\frac{\delta_1}{\delta_1 +
\varepsilon_1}
\biggr) + (\delta_2 + \varepsilon_2) \lambda\biggl(
\frac{\delta_2}{\delta_2 + \varepsilon_2} \biggr).
\]
Dividing through by $\delta_1 + \delta_2 + \varepsilon_1 + \varepsilon_2$ and
noting the relation between the arguments of $\lambda$ yields
convexity. Conversely, if $\lambda$ is strictly convex then $\kappa$ is
strictly subadditive. Showing that (\ref{ineq}) holds for all
$n>N(\delta_1,\delta_2,\varepsilon_1,\varepsilon_2)$, is equivalent to showing
%
\begin{equation}\label{ineq2}
n^{\kappa(\delta_1,\varepsilon_1)+\kappa(\delta_2,\varepsilon_2)-\kappa
(\delta
_1+\delta_2,\varepsilon_1+\varepsilon_2)} \geq L(n;\delta_1,\varepsilon_1)L(n;
\delta_2,\varepsilon_2)/L(n;\delta_1+
\delta_2,\varepsilon_1+\varepsilon_2).
\end{equation}
Taking log of each side of (\ref{ineq2}) and dividing by $\log n$ gives
$\kappa(\delta_1,\varepsilon_1)+\kappa(\delta_2,\varepsilon_2)-\kappa
(\delta_1+\delta_2,\varepsilon_1+\varepsilon_2)>0$ on the LHS and $\log\{
L(n;\delta_1,\varepsilon_1)L(n;\delta_2,\varepsilon_2)/L(n;\delta_1+\delta
_2,\varepsilon_1+\varepsilon_2)\}/\log n \to0$ as $n\to\infty$ on the RHS. Therefore,
there exists $N(\delta_1,\delta_2,\varepsilon_1,\varepsilon_2) \in
\mathbb
{N}$ such that for all $n>N(\delta_1,\delta_2,\varepsilon_1,\varepsilon_2)$
(\ref{ineq}) holds.
\end{pf}

The inequality (\ref{ineq}) has an interpretation in terms of an
existing dependence condition given in Joe \cite{Joe97}: the random vector
$(X_E,Y_E)$ gives rise to a convex $\kappa$ if asymptotically, the
conditional vector $(X_E,Y_E)|(X_E>\delta_1\log n,Y_E>\varepsilon_1\log
n)$ is \textit{more concordant} or \textit{more positively quadrant
dependent}, as $n\rightarrow\infty$, than the unconditional vector
$(X_E,Y_E)$. The majority of the examples in Section \ref{secExamples}
have a convex $\lambda$; the bivariate normal distribution with
negative dependence provides an example of a concave $\lambda$.
However, convexity is different in general to positive quadrant
dependence of $(X_E,Y_E)$, which only provides the upper bound $\lambda
(\omega)\leq1$. In a similar manner to Proposition \ref{prop3}, we
have a partial converse: $\lambda(\omega)<1$, $\omega\in(0,1)$,
implies that for each $\omega\in(0,1)$, there exists some $N(\omega)
\in\mathbb{N}$ such that for all $n>N(\omega)$,
\[
L(n;\omega,1-\omega)n^{-\lambda(\omega)} \geq n^{-\omega+
(1-\omega)} = n^{-1},
\]
that is, positive quadrant dependence for some sufficiently large
$n>N(\omega)$.

\subsubsection*{The slowly varying function}
\label{secTheSlowlyVaryingFunction}

For estimation of extreme probabilities, it is the index $\lambda
(\omega
)$ which principally controls the joint tail probability decay rate.
However, the behaviour of the slowly varying function $L$ does affect
the estimation of $\lambda$ (see Section \ref
{secStatisticalMethodology}). Owing to the assumptions of our
representation, $L$ must satisfy certain properties. One is a marginal
property: $L(n;0,\gamma)=L(n;\beta,0)=1$, by the assumption of Pareto
margins. This tells us that $L$ is identically constant iff $L \equiv
1$. A further homogeneity property arises from the assumption that
representation (\ref{A1}) holds for any $(\beta,\gamma)\in\mathbb
{R}^2_+\setminus\{\mathbf{0}\}$. For $h>0$, $L(n^h;\beta,\gamma) = L(n;
h\beta,h\gamma)$, showing we can also specify to the case $\beta
+\gamma
=1$, defining $L^*(n;\omega):= L(n;\omega,1-\omega)$.


\subsection{Generalization of scaling functions}
\label{secUnivariateRegularVariationCondition}

For any $\omega\in[0,1]$, define the transformation $T_{\omega}\dvtx
[1,\infty) \mapsto[1,\infty)^2$ by $T_{\omega}(x) = (x^\omega
,x^{1-\omega})$, and denote the joint survivor function of $(X_P,Y_P)$
by $\bar{F}\dvtx [1,\infty)^2 \mapsto(0,1]$. Assumption (\ref{A1}) can
be written
\[
\bar{F} \circ T_{\omega} \in RV_{-\lambda(\omega)} \qquad\forall \omega\in[0,1],
\]
where $RV_\rho$ stands for the class of univariate regularly varying
functions at infinity with index $\rho$ (see, e.g., Resnick \cite
{Resnick87,Resnick06}).

We can generalize the scaling functions from exact power law scaling to
regularly varying functions. Define $\bar{H}_\omega(n):=\bar{F}\circ
T_\omega(n)$, which is the survivor function of the random variable
$\min\{X_P^{1/\omega},Y_P^{1/(1-\omega)}\}$, and hence monotonic. For
$\omega=0$, we define $\min\{X_P^{1/0},Y_P^{1/1}\} =: Y_P$, and
similarly for $\omega=1$. Then we may define a non-decreasing function
$b(n)\rightarrow\infty$ such that $\bar{H}_\omega(b(n)) \sim
n^{-\lambda
(\omega)}$, $n\to\infty$, by taking $b(n) = (1/\bar{H}_\omega
)^{\leftarrow}(n^{\lambda(\omega)})$, with $(1/\bar{H}_\omega
)^{\leftarrow}$ the inverse of $1/\bar{H}_\omega$. By
Proposition 2.6(iv) and (v) of Resnick \cite{Resnick06}, we can deduce
$b(n)\in RV_1$, and thus as $n\to\infty$,
\[
\bar{F}\circ T_\omega\circ b(n) = \Prob\bigl(X_P>b(n)^{\omega
},Y_P>b(n)^{1-\omega}
\bigr) \sim n^{-\lambda(\omega)}.
\]

\subsection{Higher dimensions, strong and weak joint tail dependence}
\label{secAsymptoticDependenceAndIndependence}

The extension to higher dimensions involves a generalized notion of
asymptotic dependence and asymptotic independence. Let $\mathbf{X}$ be a
random $d$-vector with marginal distribution functions
$F_{X_1},\ldots,F_{X_d}$. We define the notion of \textit{$d$-dimensional
joint tail dependence} through the limit
%
\begin{equation}\label{jtd}
\chi(D)=\lim_{u\rightarrow1} \Prob\bigl\{F_{X_j}(X_j)>u,
\forall j \in D|F_{X_1}(X_1)>u\bigr\}
\end{equation}
for $D=\{1,\ldots,d\}$. The cases $\chi(D)>0$ and $\chi(D)=0$ define,
respectively, \textit{strong joint tail dependence} and \textit{weak joint
tail dependence} of the vector. Convergence (\ref{MSconv}) extends to
$d$-dimensions. For $\Pa$ a random $d$-vector with Pareto margins, we
have $n\Prob(\Pa/n \in\cdot)\stackrel{v}{\rightarrow} \mu(\cdot
)$ on
$[0,\infty]^d \setminus\{\mathbf{0}\}$. The homogeneous order $-1$ measure
$\mu$ has mass on the particular subcone $(0,\infty]^d \subset
[0,\infty
]^d \setminus\{\mathbf{0}\}$ only if the vector exhibits $d$-dimensional
strong joint tail dependence. For $d>2$, weak joint tail dependence
does \textit{not} imply that the limiting normalized componentwise maxima
would be mutually independent, though the converse is true. For $d=2$,
weak joint tail dependence is equivalent to asymptotic independence. In
addition, if $\chi(\{i,j\}) = 0$ for all pairs $i,j \in D$ the limiting
normalized componentwise maxima would be mutually independent.

With a slight change of notation, let $\Ex=\log\Pa$ be $d$-vectors and
$\bolds{\beta} = (\beta_1,\ldots,\beta_d)$ be the vector replacing
$(\beta,\gamma)$ in two dimensions. The function $\kappa$ is understood
to have as many arguments as the context dictates. Applying all vector
operations componentwise assume
%
\begin{equation}\label{A3}
\Prob(\Ex>\bolds{\beta}\log n) = L(n;\bolds{\beta})n^{-\kappa(\bolds{\beta})}
\end{equation}
with all quantities being defined analogously to the bivariate case.

Once more, we can define $\omega_i = \beta_i /\sum_j\beta_j$ so that
$\bolds{\omega} = (\omega_1,\ldots,\omega_d) \in S_d:=\{\mathbf{v}\in
[0,1]^d\dvtx\allowbreak\sum_{i=1}^{d}v_i = 1\}$ and the function $\lambda(\omega
_1,\ldots,\omega_{d-1}):= \kappa(\bolds{\omega})$. The previously
outlined properties of $\kappa$ and $\lambda$ hold just as in the
bivariate case; for Property \ref{property4} the upper bound follows
from extending the assumption to positive upper orthant dependence
(Joe \cite{Joe97}), other extensions are clear. In addition, we have a
consistency condition as we move across dimensions:
\[
\kappa(\beta_1,\ldots,\beta_{d})\geq\kappa(
\beta_1,\ldots,\beta_{d-1},0)= \kappa(\beta_1,\ldots,\beta_{d-1}).
\]
For a $d$-dimensional random vector, a summary of the limiting
dependence structure of the full vector is given by the values of $\chi
(C)$, for all possible $C\subseteq D$. If $\chi(C)>0$ then $\mu(B)>0$
for $B$ some Borel subset of $\bigtimes_{i=1}^d (0,\infty]^{I(i\in
C)}[0,\infty]^{I(i\notin C)}$, with $I$ the indicator function and
$\bigtimes_{i=1}^d$ representing Cartesian product of sets (e.g., if $D=\{
1,2,3\}$, $C=\{1,3\}$, the notation should be interpreted to mean
$(0,\infty]\times[0,\infty]\times(0,\infty]$). Otherwise, $\chi(C)=0$
implies $\mu(B)=0$ for all Borel $B$ on the same space, and hence the
vague convergence in (\ref{MSconv}) provides no interesting
theoretical discrimination for the dependence structure of the subset
indexed by $C$. In this sense, we have a boundary case of the limit
measure $\mu$. By contrast, $\kappa$ can still provide some
discrimination between different levels of sub-asymptotic dependence
amongst the variables indexed by~$C$.

Conversely, when $\chi(C)>0$, this strong joint tail dependence forms a
boundary case in the asymptotic theory suitable for analyzing weak
joint tail dependence. We cannot avoid this easily since we have the
dichotomy that the joint survivor function $\bar{F}$ of $\Pa$ either
decays at the same rate as the marginal survivor functions (strong
joint tail dependence) or faster than the marginal rate (weak joint
tail dependence). As a consequence here, we find that for any
distribution satisfying (\ref{A3}) and exhibiting $d$-dimensional
strong joint tail dependence, $\lambda(\omega_1,\ldots,\omega_{d-1}) =
\max_{1\leq i \leq d}\{\omega_i\}$ (see Proposition \ref{prop2}) and we
only have different forms of $\lambda$ describing different forms of
weak joint tail dependence.
%
\begin{prop}
\label{prop2} Let $\Pa$ satisfy (\ref{A3}). If $\Pa$ exhibits
$d$-dimensional strong joint tail dependence, that is, $\chi(D)>0$,
then $\kappa(\beta_1,\ldots,\beta_{d}) = \max_{1\leq i \leq
d}\{\beta_i\}$.
\end{prop}
\begin{pf}
Because $\chi(D)>0$, by definition (\ref{jtd}), we have for $\beta_1>0$
\begin{eqnarray*}
&&
\lim_{n\rightarrow\infty}\Prob(X_{P,2}>n,\ldots,X_{P,d}>n|X_{P,1}>n) \\
&&\quad=
\lim_{n\rightarrow\infty}\Prob\bigl(X_{P,2}>n^{\beta_1},\ldots,X_{P,d}>n^{\beta_1}|X_{P,1}>n^{\beta_1}\bigr)
= \chi(D) > 0.
\end{eqnarray*}
Without loss of generality, assume $\beta_1=\max_{1\leq i \leq d}\{
\beta_i\}$. Then we also have
\begin{eqnarray*}
1 &\geq&\Prob\bigl(X_{P,2}>n^{\beta_{2}},\ldots,X_{P,d}>n^{\beta
_{d}}|X_{P,1}>n^{\beta_1} \bigr)\\
&\geq&\Prob\bigl(X_{P,2}>n^{\beta_1},\ldots,X_{P,d}>n^{\beta
_1}|X_{P,1}>n^{\beta_1}\bigr)
\rightarrow\chi(D)
\end{eqnarray*}
as $n\rightarrow\infty$. Hence, $\Prob(X_{P,1}>n^{\beta
_1},X_{P,2}>n^{\beta_{2}},\ldots,X_{P,d}>n^{\beta_{d}}) \asymp
n^{-\beta
_1}$, as $n\to\infty$, that is, $\kappa(\beta_1,\ldots,\beta_{d})
= \max_{1\leq i \leq d}\{\beta_i\}$.
\end{pf}

In practice as $d$ grows it becomes more likely that a scenario will
arise for which some sets $C \subseteq D$ yield $\chi(C)=0$ and some
yield $\chi(C)>0$, and as such, tools from both dependence paradigms
will be useful. For relative simplicity of exposition, we return the
focus to the bivariate case for the rest of the paper.


\section{Bivariate limits}
\label{secConditionalMultivariateLimits} In this section, we detail
some results concerning bivariate limits, initially by conditioning
upon the event that $(X_P,Y_P)
\in(n^\omega,\infty)\times(n^{1-\omega},\infty)$.

\subsection{Assumption}
\label{secAssumption}

Consider, for $x,y\geq1$, $\Prob(X_P>n^{\omega} x,Y_P>n^{1-\omega}
y |
X_P>n^{\omega},Y_P>n^{1-\omega})$; under assumption (\ref{A1}) this is
given by
%
\begin{eqnarray}\label{cond1}
&&
\frac{L (n; \omega+ {\log x}/{\log n},1-\omega+
{\log
y}/{\log n} )}{L(n;\omega,1-\omega)}\nonumber\\[-8pt]\\[-8pt]
&&\quad{}\times\exp\biggl[\log n \biggl\{ \kappa
(\omega,1-\omega) - \kappa
\biggl(\omega+ \frac{\log x}{\log
n},1-\omega+ \frac{\log y}{ \log n} \biggr) \biggr\}
\biggr].\nonumber
\end{eqnarray}
We shall initially explore the limiting behaviour of this expression
under the assumption that $\kappa$ is differentiable and that
%
\begin{equation}\label{A2}
\lim_{n\rightarrow\infty} \frac{L (n; \omega+{\log
x}/{\log
n},1-\omega+{\log y}/{\log n} )}{L(n; \omega, 1-\omega)} = 1.
\end{equation}
Such assumptions are smoothness conditions on the behaviour of the
joint distributional tail over different $\omega$. They are analogous
to supposing ray independence of the Ledford and Tawn \cite
{LedfordTawn97} paradigm; a
condition which actually holds very widely including for all of the
asymptotically independent examples of the Ledford and Tawn \cite
{LedfordTawn97} paper.
The development under these assumptions will permit us to observe how a
generalized notion of multivariate regular variation is possible at the
end of the section.


\subsection{Conditional limits under the assumption}
\label{secConditionalLimitsUnderTheAssumption}

Since we assume that $\kappa$ is differentiable at the point $\bolds{\omega
}:=(\omega,1-\omega)$, then
\[
\lim_{n\rightarrow\infty} \log n \biggl[\kappa(\omega,1-\omega) -
\kappa\biggl(
\omega+ \frac{\log x}{\log n},1-\omega+ \frac
{\log
y}{\log n} \biggr) \biggr]= -\nabla
\kappa(\omega,1-\omega)\cdot(\log x, \log y)^T
\]
with $\nabla\kappa$ the gradient vector of $\kappa$. Therefore, the
conditional probability limit from equation (\ref{cond1}) is
%
\begin{equation}\label{NHPar}
\exp\bigl\{-\nabla\kappa(\omega,1-\omega)\cdot(\log x, \log y)^T
\bigr\} = x^{-\kappa_1(\omega)}y^{-\kappa_2(\omega)},
\end{equation}
that is, independent Pareto variables with shape parameters
\[
\bigl\{\kappa_1(\omega),\kappa_2(\omega)\bigr\}:=
\biggl(\frac
{\partial\kappa}{\partial\beta},\frac{\partial\kappa}{\partial
\gamma
} \biggr)\bigg|_{(\omega,1-\omega)} = \bigl\{
\lambda(\omega)+(1-\omega)\lambda'(\omega), \lambda(\omega) -
\omega\lambda'(\omega) \bigr\}.
\]
By Property \ref{property2}, it is clear that $\kappa_1(\omega
),\kappa_2(\omega)\geq0$. In exponential margins, for $x,y\geq0$, this reads
%
\begin{eqnarray}\label{NHExp}
\lim_{n\rightarrow\infty}\frac{\Prob\{X_E-\omega\log
n>x,Y_E-(1-\omega
)\log n>y\}}{\Prob\{X_E>\omega\log n,Y_E>(1-\omega)\log n\}} &=& \exp
\bigl\{ -\nabla\kappa(\omega,1-
\omega) \cdot(x, y)\bigr\} \quad\nonumber\\[-8pt]\\[-8pt]
&=& \exp\bigl\{-\kappa_1(\omega) x-
\kappa_2(\omega) y\bigr\}.\nonumber
\end{eqnarray}
Although stochastic independence of the conditioned normalized
variables arises in the limit, limits (\ref{NHPar}) and (\ref{NHExp})
are interpretable in terms of the original degree of dependence at the
finite level. This can be observed by the relation $\kappa_1(\omega) =
\kappa_2(\omega)+\lambda'(\omega)$. As the variables become
independent, $\lambda(\omega)$ and $\lambda'(\omega)$ approach 1
and 0,
respectively, thus $\kappa_1(\omega)$ and $\kappa_2(\omega)$ approach
equality. As dependence increases, then $\kappa_1(\omega)$ decreases
relative to $\kappa_2(\omega)$ for $\omega<1/2$, and vice versa.

The boundary case of asymptotic dependence has $\kappa_1(\omega)=0$ and
$\lambda'(\omega) = -1$ for $\omega<1/2$ and $\kappa_2(\omega)=0$ and
$\lambda'(\omega) = 1$ for $\omega>1/2$. Under asymptotic dependence,
the function $\lambda(\omega)=\max\{\omega,1-\omega\}$ is not
differentiable at $\omega=1/2$, and the above limits (\ref{NHPar})
and (\ref{NHExp}) do not hold on this line. However, when $\omega\neq
1/2$, we can still get the conditional limit from (\ref{NHPar}) to be
$x^{-1} I(\omega>1/2)+y^{-1} I(\omega<1/2)$. When $\omega=1/2$,
standard theory for asymptotically dependent distributions provides
that the limit is $\{x^{-1}+y^{-1} - V(x,y)\}/\{2-V(1,1)\}$, where
$V(x,y):=\mu(\{[0,x]\times[0,y]\}^c)$ for the appropriate limiting
measure $\mu$ in convergence (\ref{MSconv}).

The asymptotic dependence case highlights some subtleties that occur
when the assumptions of differentiable $\kappa$ and assumption (\ref
{A2}) break down. In this case, the function $L$ inherits additional
structure at the point where $\kappa$ is not differentiable, and
classical multivariate regular variation conditions become the
appropriate tools for analysis.

\subsection{Multivariate regular variation}
\label{secMultivariateRegularVariation}

The weak convergence of probability measures given by limit (\ref
{NHPar}) can be expressed more generally as
%
\begin{equation}\label{Condconv}
\Prob\biggl\{ \biggl(\frac{X_P}{n^\omega},\frac{Y_P}{n^{1-\omega
}} \biggr) \in\cdot\bigg|
\biggl(\frac{X_P}{n^\omega},\frac{Y_P}{n^{1-\omega
}} \biggr)\in[1,\infty)^2
\biggr\} \stackrel{w} {\rightarrow} \pi (\cdot;\omega)
\end{equation}
on $(1,\infty)^2$ as $n\rightarrow\infty$, with the limit measure
$\pi
(\cdot;\omega)$ homogeneous of order $-\{\kappa_1(\omega)+\kappa
_2(\omega)\}$. This reveals an alternative non-standard type of regular
variation condition. A function $q$ is (standard) multivariate
regularly\vspace*{1pt} varying with index $\rho\in\mathbb{R}$ and limit function
$\Lambda$ if $\lim_{t\rightarrow\infty} q(t\mathbf{x})/q(t\mathbf{1}) =
\Lambda
(\mathbf{x})$, $\mathbf{x}\in\mathbb{R}^d_+$, the function $q(t\mathbf{1})
=:Q(t)\in RV_\rho$, is univariate regularly varying, and $\Lambda
(c\mathbf
{x}) = c^\rho\Lambda(\mathbf{x})$, $c>0$. The common scaling $t\in RV_1$
can be replaced by any common scaling $b(t)\in RV_{\alpha}$, $\alpha
>0$, yielding an unchanged limit function and $Q\circ b(t)\in RV_{\rho
\alpha}$. In limit (\ref{Condconv}) however, different scalings are
applied to each component, in the classes $RV_\omega$ and
$RV_{1-\omega
}$, $Q\circ T_\omega(t) \in RV_{-\lambda(\omega)}$, and the order of
homogeneity of the limit measure is $-\{\kappa_1(\omega)+\kappa_2(\omega
)\}$. When $\omega=1/2$, the common scaling returns standard
multivariate regular variation since the limit measure is homogeneous
of order $-\{\kappa_1(1/2)+\kappa_2(1/2)\}=-\kappa(1,1)=-2\lambda
(1/2)=-1/\eta$. Note this is a different form of non-standard regular
variation to that given in Resnick \cite{Resnick06}, as there the non-equal
scalings are due to non-equal margins.

Limit (\ref{Condconv}) provides an asymptotic link between the
probabilities of lying in Borel sets $B\subset[1,\infty)^2$ where each
component is scaled by $(n^\omega,n^{1-\omega})$ as
\[
\Prob\biggl\{ \biggl(\frac{X_P}{n^\omega},\frac{Y_P}{n^{1-\omega
}} \biggr) \in
\bigl(t^\omega,t^{1-\omega}\bigr)B \biggr\} \sim t^{-\lambda(\omega
)}\Prob
\biggl\{ \biggl(\frac{X_P}{n^\omega},\frac{Y_P}{n^{1-\omega}} \biggr)
\in B \biggr\}
\]
as $n\rightarrow\infty$, $t>1$, with vector arithmetic applied
componentwise. This follows since $\omega\kappa_1(\omega) +
(1-\omega
)\kappa_2(\omega) = \kappa(\omega,1-\omega) = \lambda(\omega)$ by
homogeneity. Together with the corresponding relation in exponential
margins, this provides the analogous equations to (\ref{LText}).

By the form of limit (\ref{Condconv}), we observe that this can be
modified to standard multivariate regular variation. Define the
bijective map $U_\omega\dvtx (1,\infty)^2 \mapsto(1,\infty)^2$ by
$U_\omega(x,y) = (x^{\omega}, y^{1-\omega})$, for any fixed $\omega
\in(0,1)$. Then $\bar{F}\circ U_{\omega}$ is multivariate regularly
varying of index $-\lambda(\omega)$. Equivalently we can think of this
as standard multivariate regular variation of the law of the random
vector $U_\omega^{\leftarrow}(X_P,Y_P)=(X_P^{1/\omega},
Y_P^{1/(1-\omega
)})$, $\omega\in(0,1)$. In place of convergence (\ref{Condconv}), we
can thus consider,
%
\begin{equation}\label{Condconv2}
\Prob\bigl\{U_\omega^{\leftarrow}(X_P,Y_P)/n
\in\cdot| U_\omega^{\leftarrow}(X_P,Y_P) / n
\in[1,\infty)^2 \bigr\} \stackrel{w} {\rightarrow} \pi ^*(\cdot;
\omega)
\end{equation}
on $(1,\infty)^2$, where for Borel $B\subset(1,\infty)^2$, $c>0$,
$\pi ^*(c B;\omega) = c^{-\lambda(\omega)}\pi ^*(B;\omega)$. When
$\omega
=1/2$, this is standard multivariate regular variation, and does not
therefore require the assumptions of differentiable $\kappa$
and (\ref{A2}). Consequently, (\ref{Condconv2}) provides a general
multivariate
asymptotic characterization for both asymptotically dependent and
asymptotically independent random vectors, where the structure rests in
$\pi ^*(\cdot;1/2)$ under asymptotic dependence, and the combination of
$\lambda(\omega)$ and $\pi ^*(\cdot;\omega)$ under asymptotic
independence.

As in the univariate case, the scaling functions can be generalized
from exact powers of $n$ to regularly varying functions. Under our
current assumptions, we can find $a(n)\in RV_{1/\lambda(\omega)}$ such
that for Borel $B\subset(1,\infty)^2$,
\[
\lim_{n\to\infty} n\Prob\bigl\{U_\omega^{\leftarrow}(X_P,Y_P)
/ a(n) \in B \bigr\} = \pi ^*(B; \omega)
\]
with $\pi ^*(\partial B)=0$ and $\pi ^*$ as given in convergence (\ref
{Condconv2}). This follows by defining $a(n) = (1/\bar{H}_\omega
)^{\leftarrow}(n)$, as in Section \ref
{secUnivariateRegularVariationCondition}, and using the assumptions of
Section \ref{secAssumption} on $L$ and $\kappa$, or standard
multivariate regular variation assumptions.


\section{Connections to existing theory for asymptotic independence}
\label{secConnectionsToExistingTheory}

In this section, we detail more carefully the connections between our
representation (\ref{A1}) and the existing theories which permit
non-trivial treatment of asymptotic independence. Specifically, these
are the theories introduced by Ledford and Tawn \cite{LedfordTawn97} and
Heffernan and Tawn \cite{HeffernanTawn04}. These representations were
originally phrased in
terms of standard Fr\'{e}chet and standard Gumbel margins,
respectively. Here, we continue to use the asymptotically equivalent
standard Pareto and exponential margins.


\subsection{Coefficient of tail dependence}
\label{secCoefficientOfTailDependence}
The assumption of Ledford and Tawn \cite{LedfordTawn97} is given
by (\ref{LT}).
Supposing that $\kappa$ is differentiable then the results of
Section \ref{secConditionalMultivariateLimits} allow us to write
%
\begin{eqnarray}\label{LTlink}
&&
\Prob\bigl(X_P>n^{\omega} x,Y_P>n^{1-\omega}
y\bigr) \nonumber\\[-8pt]\\[-8pt]
&&\quad= n^{-\lambda(\omega
)}x^{-\kappa_1(\omega)}y^{-\kappa_2(\omega)}L(n;\omega+ \log x/\log
n,1-\omega+ \log y/\log n)\nonumber
\end{eqnarray}
for $x,y \geq1$ as $n\rightarrow\infty$. Setting $\omega=1/2$ in
equation (\ref{LTlink}) returns the set-up of equation (\ref{LT}). From
this, we see that the representations are equivalent when $c_1 = \kappa
_1(1/2) = \lambda(1/2) + \lambda'(1/2)/2$, $c_2 = \kappa_2(1/2) =
\lambda(1/2)-\lambda'(1/2)/2$ and $\Lsv(nx,ny) \sim L(n;1+ \log
x/\log
n,1+ \log y/\log n)$, $n\to\infty$.

Ramos and Ledford \cite{RamosLedford09} introduced an expanded
characterization of the
Ledford and Tawn \cite{LedfordTawn97} framework. Similarly to
Section \ref
{secConditionalMultivariateLimits}, they consider convergence of
measures $\Prob(X_F>nx,Y_F>ny | X_F>n,Y_F>n)$, with $(X_F,Y_F)$
marginally standard Fr\'{e}chet, however without the condition (\ref
{A2}). This entails the possibility of more structure in the limit
measure, and is a more natural idea under equal marginal growth rates,
where any non-constant limit in equation (\ref{A2}) would be
homogeneous of order 0. Equation (\ref{Condconv2}) has the
Ramos and Ledford \cite{RamosLedford09} assumption as a special case.
However, their focus is
on developing consequences of the multivariate regular variation
condition under common scaling of the margins, whilst our focus remains
on consequences of generalizing the scaling.


\subsection{Conditioned limit theory}
\label{secConditionedLimitTheory}
Heffernan and Tawn \cite{HeffernanTawn04} generated a limit theory for
the distribution
of a variable $X$, suitably normalized, conditional upon the
concomitant variable, $Y$, being extreme. In standard Gumbel margins,
they considered non-degenerate limits of
%
\begin{equation}\label{cevht}
\Prob\bigl(\bigl\{X_G-b(y)\bigr\}/a(y)\leq x | Y_G =
y \bigr)
\end{equation}
as $y\rightarrow\infty$. The normalizing functions $a,b$ satisfy
$a(y)>0$, $b(y) \leq y$, and give unique limit distributions up to type
(see Heffernan and Tawn \cite{HeffernanTawn04}, Theorem 1). Precise
forms of $a,b$ are
stated on the understanding that modifications which do not change the
type of limit distribution are allowable. In this limit theory,
asymptotic dependence is once more a boundary case requiring the
normalizations $a(y) = 1$ and $b(y)=y$. Formulation (\ref{cevht}) was
elaborated on by Heffernan and Resnick \cite{HeffernanResnick07}, who
described limit theory
for both
%
\begin{eqnarray}\label{cevhr}
&&\Prob\bigl(\bigl\{X_E-b(\log n)\bigr\}/a(\log n)\leq x|
Y_E > \log n \bigr)
\nonumber
\\[-8pt]
\\[-8pt]
\nonumber
&&\Prob\bigl(\bigl\{X_E-b(Y_E)
\bigr\}/a(Y_E)\leq x| Y_E > \log n \bigr)
\end{eqnarray}
as $n \rightarrow\infty$. The formulation of Heffernan and Resnick
\cite{HeffernanResnick07}
is in fact slightly more general as the marginal distribution of the
$X$ variable is left unspecified, however it is convenient here to
suppose a particular form.

It is possible to draw some modest links between our representation and
the limits in (\ref{cevhr}). However, the full structure of the
function $L$ in representation (\ref{A1}) across different $\beta,\gamma
$ is key in determining the limits of (\ref{cevhr}). By fixing $(\beta
,\gamma)\in\mathbb{R}_+^2\setminus\{\mathbf{0}\}$, and considering
asymptotics in $n$, then any slowly varying function $L_1(n;\beta
,\gamma
)\sim L(n;\beta,\gamma)$, $n\to\infty$, suffices for our calculations:
this highlights the difference between the theories.

Under assumption (\ref{A1}) and with $\gamma=1$,
%
\begin{eqnarray}
\label{cond} \Prob(X_E>\beta\log n|Y_E> \log n) &=& L(n;
\beta,1) n^{1-\kappa
(\beta,1)}\nonumber\\[-8pt]\\[-8pt]
&=& L(n; \beta,1)\exp\bigl[\bigl\{1-\kappa(\beta,1)\bigr\}
\log n\bigr]\nonumber
\end{eqnarray}
as $n\rightarrow\infty$. If we view the $\log n$ in the first argument
of the conditional probability (\ref{cond}) as being part of the
conditioning variable, then we want to find $\beta= \beta(n,x)$ such
that $\beta(n,x)\log n = a(\log n)x+b(\log n)$ gives a non-degenerate
limit in $x$ of equation (\ref{cond}) as $n \rightarrow\infty$. That
is, we are searching for the functions $a,b,k$ such that
%
\begin{equation}\label{condmarg}
\lim_{n\rightarrow\infty}\Prob\bigl(\bigl\{X_E-b(\log n)\bigr\} /a(\log
n)>x|Y_E> \log n \bigr) = \exp\bigl\{-k(x)\bigr\},
\end{equation}
defines a valid survivor function. Supposing that $\lim_{n\rightarrow
\infty}\beta(n,x) = 0$, one can observe that the only simple link
between the two theories arises when $L\equiv1$, since by the marginal
condition $L(n;\beta,0)=L(n;0,\gamma)=1$. In this case,
equation (\ref{cond}) shows that the components required for the
limit (\ref{condmarg}) satisfy the relation that
%
\begin{equation}\label{beta}
\kappa\bigl(\beta(n,x),1\bigr) - 1 \sim k(x)/\log n
\end{equation}
as $n\rightarrow\infty$. When $L \not\equiv1$, then the first order
theory for $\Prob(X_E>\beta\log n, Y_E> \gamma\log n)$ that is the
focus of our work is insufficient to yield the full limit distribution
of equation (\ref{condmarg}). Example \ref{secMorgenstern} of the
following section provides an illustration of this.


\subsection{Examples}
\label{secExamples}

Below we provide examples to illustrate the theory developed and the
links with existing representations. In each case, we provide the
copula, $C(u,v)$, of the distribution (i.e., the joint distribution
function on standard uniform margins) and the asymptotic joint survivor
function (\ref{A1}) on exponential margins. Table \ref{tabexamples}
summarizes the previously defined quantities $\kappa(\beta,\gamma)$;
%
\begin{table}
\caption{Summary of important theoretical quantities for
representations of Examples
\protect\ref{secBivariateNormalDistribution}--\protect\ref{secLowerJointTailOfTheClaytonDistribution}
under the theory of the current paper, Ledford and Tawn
\cite{LedfordTawn96} and Heffernan and Tawn
\cite{HeffernanTawn04}/Heffernan and Resnick \cite{HeffernanResnick07}}
\label{tabexamples}
\begin{tabular*}{\tablewidth}{@{\extracolsep{\fill}}llllll@{}}
\hline
Example &$\kappa(\beta,\gamma)$ & $\eta$ & $a(\log n)$ & $b(\log n)$ &
$\exp\{-k(x)\}$ \\
\hline
\ref{secBivariateNormalDistribution}& (a) $\frac{\beta+\gamma- 2\rho
(\beta\gamma)^{1/2}}{1-\rho^2}$; & $(1+\rho)/2$ & $(2\rho^2\log
n)^{1/2}$\tabnoteref{tb}
& $\rho^2 \log n$\tabnoteref{tb} &
$1-\Phi\bigl\{\frac{x}{(1-\rho^2)^{1/2}} \bigr\}$\\[2pt]
& (b) $\max\{\beta,\gamma\}$\\[2pt]
\ref{secInvertedBivariateExtremeValueDistribution} & $V(1/\beta,
1/\gamma)$ & $1/V(1,1)$ & \tabnoteref{ta} & \tabnoteref{ta} & \tabnoteref{ta}\\[2pt]
\ref{secMorgenstern} & $\beta+\gamma$ & $1/2$ & $1$ & $0$ &
$e^{-x}[1+\alpha(1-e^{-x})]$\\[2pt]
\ref{secBivariateExtremeValueDistribution} & $\max\{\beta,\gamma\}$ &
$1$ & $1$ & $\log n$ & $e^{-x}[1+e^{x}-V(1,e^{-x})]$\\[2pt]
\ref{secLowerJointTailOfTheClaytonDistribution} & $\max\{\beta
,\gamma
\}$ & $1$ & $1$ & $\log n$ & $[e^{x/\alpha} +1]^{-\alpha}$\\
\hline
\end{tabular*}
\tabnotetext[\mbox{$*$}]{ta}{Denotes
that the explicit form depends upon $V$.}
\tabnotetext[\mbox{$\dagger$}]{tb}{Denotes for
$\rho>0$. Example \ref{secBivariateNormalDistribution}(a): $\rho^2<\min \{\beta/\gamma,\gamma/\beta\}$ or $\rho<0$
and $\min \{\beta,\gamma\}>0$; Example \ref{secBivariateNormalDistribution}(b): all other $\rho$.}
\end{table}
$\eta$; $a(\log n), b(\log n)$; and $\exp\{-k(x)\}$. Detailed
derivations are provided for two interesting examples: the bivariate
normal, and the inverted bivariate extreme value distribution. Other
bivariate examples are stated more briefly. We additionally present a
trivariate example for which $\chi(\{1,2,3\})=0$ but some $\chi(\{
i,j\}
)>0$ to illustrate behaviour both across dimensions and extremal
dependence classes.

\begin{example}[(Bivariate normal distribution)]
\label{secBivariateNormalDistribution}
\[
C(u,v) = \Phi_2\bigl(\Phi^{-1}(u),\Phi^{-1}(v);
\rho\bigr),\qquad \rho\in[-1,1],
\]
with $\Phi$ the standard normal distribution function and $\Phi_2(\cdot
,\cdot;\rho)$ the standard bivariate normal distribution function.
Denote also by $\phi_2$ the associated joint density. The survivor
function in this case is defined only via an integral, however an
asymptotic expansion yields $\kappa$ and the leading order term of $L$.
Let $(X_N,Y_N)$ denote the vector on standard normal margins. Interest
lies in $\Prob(X_E>\beta\log n,Y_E>\gamma\log n) = \Prob\{X_N>\Phi^{-1}
(1-n^{-\beta}), Y_N > \Phi^{-1}(1-n^{-\gamma})\}$. Let $x_\beta:=
\Phi^{-1}
(1-n^{-\beta})$; $x_\gamma:= \Phi^{-1}(1-n^{-\gamma})$. Then as
$n\to
\infty$
\[
x_\beta= (2\beta\log n)^{1/2} -\bigl(\log\log
n^\beta+\log4\uppi \bigr)/\bigl\{ {2(2\beta\log n)^{1/2}}\bigr\}
+\mathrm{o} \bigl\{(\log n)^{-1/2} \bigr\}
\]
and similarly $x_\gamma$, thus $x_\beta\sim(2\beta\log n)^{1/2}$,
$x_\gamma\sim(2\gamma\log n)^{1/2}$, and hence $x_\beta\sim(\beta
/\gamma)^{1/2} x_\gamma$ as $n\rightarrow\infty$ for $\beta,\gamma>0$.
Bounds on the multivariate Mills' ratio $ \Prob(X_N>x_\beta
,Y_N>x_\gamma
) / \phi_2(x_\beta,\allowbreak x_\gamma; \rho)$, obtained by Savage \cite{Savage62}, or
the asymptotic results of Ruben \cite{Ruben64}, see also
Hashorva and H{\"u}sler~\cite{HashorvaHusler02}, provide that for
$x_\beta>\rho x_\gamma$ and
$x_\gamma>\rho x_\beta$,
\[
\Prob(X_N>x_\beta,Y_N>x_\gamma) \sim
L_1(n;\beta,\gamma) n^{-
\{
({\beta+\gamma-2\rho(\beta\gamma)^{1/2}})/({1-\rho^2}) \}}
\]
as $n\rightarrow\infty$ with
%
\begin{eqnarray}\label{nlsv}
L_1(n;\beta,\gamma) &=& (4\uppi \log n )^{({2
\rho^2 - \rho\{ (\beta/\gamma)^{1/2}+(\gamma/
\beta)^{1/2}\}})/({2(1-\rho^2)})}\nonumber\\[-8pt]\\[-8pt]
&&{}\times\frac
{\beta^{({1-\rho(\gamma/\beta)^{1/2}})/({2(1-\rho^2)})}\gamma^
{({1-\rho(\beta
/\gamma)^{1/2}})/({2(1-\rho^2)})}(1-\rho^2)^{3/2}}{(\beta^{1/2}-\rho
\gamma^{1/2})(\gamma^{1/2}-\rho\beta^{1/2})}.\nonumber
\end{eqnarray}
Note that $L_1\neq L$, that is, higher order terms of $L$ have been
neglected, but $L_1/L\to1$ as $n\to\infty$. In terms of $\beta
,\gamma
$, if $\rho>0$, $\rho^2<\min\{\beta/\gamma,\gamma/\beta\}$
implies $\rho
<\min\{x_\beta/x_\gamma,x_\gamma/x_\beta\}$ for all sufficiently large
$n$, so that the above results hold. For $\rho<0$, the above holds for
any $\beta>0$, \mbox{$\gamma>0$}. The necessity of the strict inequalities arises
from the assumption that $x_\beta$ and $x_\gamma$ are both growing;
indeed one may observe that setting, for example, $\beta=0$, with $\rho<0$,
does not yield the appropriate marginal distribution for this
representation, and more careful analysis is required. When $\rho>0$
and $\beta\leq\rho^2\gamma$ or $\gamma\leq\rho^2\beta$, it is
simple to
derive the results directly. We have
%
\begin{eqnarray}\label{jtprobnorm}
&&\Prob(X_N>x_\beta, Y_N > x_\gamma)\nonumber\\
&&\quad= \int_{x_\gamma}^{\infty} \biggl\{ 1-\Phi\biggl(
\frac{x_\beta- \rho y}{(1-\rho^2)^{1/2}} \biggr) \biggr\} \frac{1}{(2\uppi
)^{1/2}}\exp\bigl
\{-y^2/2\bigr\} \mrmd y
\\
&&\quad= \frac{\phi(x_\gamma)}{x_\gamma}\int_{0}^{\infty} \biggl\{ 1-
\Phi\biggl(\frac{x_\beta-\rho x_\gamma-\rho z/x_\gamma}{(1-\rho
^2)^{1/2}} \biggr) \biggr\}\exp\bigl\{-z-z^2/2x_\gamma^2
\bigr\} \mrmd z,\nonumber
\end{eqnarray}
the second line following with a change of variables $z=x_\gamma
(y-x_\gamma)$. It is clear that $\beta< \rho^2\gamma$ ensures that the
argument of $\Phi$ tends to $-\infty$, so that by dominated
convergence, the integral converges to $\int_0^\infty e^{-z}\mrmd z = 1$.
If $\beta= \rho^2\gamma$, then $x_\beta-\rho x_\gamma= (\rho-\rho
^{-1})\log\log n / \{2(2\gamma\log n)^{1/2}\} + \mathrm{O}\{(\log n)^{-1/2}\}$
as $n\rightarrow\infty$. Consequently, the argument of $\Phi$ tends to
0, and the integral converges to $\int_0^\infty e^{-z}/2\mrmd z = 1/2$.
Overall therefore $\Prob(X_N>x_\beta,Y_N>x_\gamma) \sim L_1(n;\beta
,\gamma) n^{-\kappa(\beta,\gamma)}$ with
%
\begin{eqnarray}
\label{L1normal} \kappa(\beta,\gamma) &=& \cases{ \displaystyle \frac{\beta+\gamma
-2\rho(\beta\gamma)^{1/2}}{1-\rho^2}, &\quad $
\displaystyle \rho^2 < \min\{\beta/\gamma,\gamma/\beta\}$ or $\rho<0$\cr
&\quad and $\min\{
\beta,\gamma\}>0$,
\cr
\displaystyle \max\{\beta,\gamma\}, &\quad otherwise,}
\nonumber\\[-8pt]\\[-8pt]
L_1(n;\beta,\gamma) &=& \cases{ \mbox{(\ref{nlsv})},&\quad $
\displaystyle \rho^2 < \min\{\beta/\gamma,\gamma/\beta\}$ or $\rho<0$\vspace*{1pt}\cr
&\quad and $\min\{
\beta,\gamma\}>0$,
\vspace*{1pt}\cr
1, &\quad $\displaystyle \rho^2 > \gamma/\beta$ or $
\rho^2 > \beta/\gamma$,
\cr
1/2, &\quad $\displaystyle \rho^2 = \beta/\gamma$
or $\rho^2 = \gamma/\beta$.}
\nonumber
\end{eqnarray}
The function $L_1$ is such that $L(n;\beta,\gamma) \sim L_1(n;\beta
,\gamma)$ in each case, and whilst the form of $L_1$ may be written as
discontinuous, higher order terms in $L$ ensure that the joint survivor
function behaves smoothly across the quadrant.
The function $\kappa$ is both continuous and differentiable on the
lines $\rho^2=\beta/\gamma$, $\rho^2=\gamma/\beta$, and attains the
expected limits $\max\{\beta,\gamma\}$ and $\beta+\gamma$ as $\rho
\rightarrow1$ and $\rho\rightarrow0$, respectively.

The normalizing functions $a,b$ and limit distribution for the
Heffernan and Tawn \cite{HeffernanTawn04} representation can be found
through consideration of
equation (\ref{jtprobnorm}) in the case where $\rho>0$. For $\rho=0$
the variables are independent thus the required normalization is
trivially $a(\log n)= 1$, \mbox{$b(\log n)= 0$}. For $\rho<0$ the expansion used
to find $\kappa$ is not sufficient, since the variables were assumed to
be growing in both margins; in fact to get a non-degenerate conditional
limit under negative dependence, one needs to consider the case where
the margins are not growing simultaneously. In exponential margins,
equation (\ref{jtprobnorm}) yields
\begin{eqnarray*}
&&
\Prob\bigl(X_E>\beta(n,x)\log n| Y_E > \log n\bigr)\\
&&\quad
\sim\int_{0}^{\infty
} \biggl\{ 1-\Phi\biggl(
\frac{x_{\beta(n,x)} -\rho x_1 -\rho z/x_1}{(1-\rho^2)^{1/2}} \biggr)
\biggr\}\exp\bigl\{-z-z^2/2x_1^2
\bigr\} \mrmd z.
\end{eqnarray*}
Non-degeneracy of this in the limit requires that $\beta(n,x)^{1/2}
-\rho\asymp(\log n)^{-1/2}, n\rightarrow\infty$. The precise
normalization required is unique up to type only. Recalling $\beta(n,x)
= \{a(\log n)x + b(\log n)\}/\log n$, one possibility is to solve
\[
\biggl\{\frac{a(\log n)}{\log n}x + \frac{b(\log n)}{\log n} \biggr\}
^{1/2} -\rho
\sim\frac{x}{(2\log n)^{1/2}}
\]
as $n\rightarrow\infty$. This provides $b(\log n)=\rho^2\log n + \mathrm{o}\{
(\log n)^{1/2}\}$, and $a(\log n) = \rho(2\log n)^{1/2}+ \mathrm{o}\{(\log
n)^{1/2}\}$. Taking\vspace*{1pt} this choice of $a, b$ leads to the limiting
survivor function $1-\Phi\{x/(1-\rho^2)^{1/2}\}$. This normalization
and limit is the same as that stated in Heffernan and Resnick~\cite
{HeffernanResnick07}; in
Heffernan and Tawn \cite{HeffernanTawn04} the normalization $a(y) =
y^{1/2}$, $b(y)=\rho^2 y$ leads to the limiting survivor function
$1-\Phi\{x/(2\rho^2(1-\rho^2))^{1/2}\}$. Note that taking specifically
$a(\log n) = \rho(2\log
n)^{1/2} + 1 $ avoids the concerns raised in Heffernan and Tawn \cite
{HeffernanTawn04} of
discontinuity in normalization as $\rho\searrow0$.
\end{example}
%
\begin{example}[(Inverted bivariate extreme value distribution)]
\label{secInvertedBivariateExtremeValueDistribution}
\[
C(u,v) = u + v - 1 + \exp\bigl[-V\bigl\{-1/\log(1-u),-1/\log(1-v)\bigr\}
\bigr],
\]
where $V$ is as defined in Section \ref
{secConditionalLimitsUnderTheAssumption}, and termed the \textit
{exponent function} of a bivariate extreme value distribution. This gives
\[
\Prob(X_E>\beta\log n,Y_E>\gamma\log n) = \exp\bigl
\{-V(1/\beta\log n, 1/\gamma\log n)\bigr\} = n^{-V(1/\beta, 1/\gamma)}.
\]

General forms of the normalization functions and limit distribution
cannot be expressed explicitly without knowledge of $V$. However, one
can consider general classes of exponent functions $V$ and characterize
the limits for these classes. The function $V$ can be expressed as
\[
V(x,y) = \int_{0}^{1} \max\bigl\{w / x,
(1-w)/y \bigr\}H(\mrmdd w),
\]
where $H$ is a measure satisfying $\int_{0}^{1} w H(\mrmdd w) = \int_{0}^{1}
(1-w) H(\mrmdd w) = 1$ (e.g., Beirlant \textit{et al.}~\cite{Beirlant04}).
Suppose $H$ has density
$h(w)$ with respect to Lebesgue measure on $(0,1)$, such that $h(w)
\sim s w^t$ as $w\rightarrow0$, $h(w) \sim u (1-w)^v$ as $w\rightarrow
1$, for $s,u>0$, $t,v>-1$, with no mass at $w=0, w=1$. Then
\[
V(1/\beta, 1) = \beta\int_{(1+\beta)^{-1}}^1 w H(\mrmdd w) + 1-
\int_{(1+\beta)^{-1}}^1 (1-w) H(\mrmdd w).
\]
As $\beta\rightarrow0$, this can be written
\begin{eqnarray*}
V(1/\beta, 1) -1 &\sim& \beta u \int_{(1+\beta)^{-1}}^1 w
(1-w)^v \mrmd w - u\int_{(1+\beta)^{-1}}^1
(1-w)^{v+1} \mrmd w
\\
& = &\frac{u}{(v+1)(v+2)}\beta^{v+2} + \RMo\bigl(\beta^{v+2}\bigr).
\end{eqnarray*}
Since $L\equiv1$ we can now use equation (\ref{beta}) to obtain the
normalizations and the limit:
\[
\frac{u}{(v+1)(v+2)}\beta(n,x)^{v+2} \sim V\bigl(\beta(n,x)^{-1},
1\bigr) -1 = \kappa\bigl(\beta(n,x),1\bigr) -1 \sim k(x)/\log n,
\]
so that $a(\log n), b(\log n)$ and $k(x)$ solve
\[
\frac{u}{(v+1)(v+2)} \biggl\{\frac{a(\log n)x + b(\log n)}{\log
n} \biggr\}^{v+2} \sim k(x)/
\log n
\]
as $n\rightarrow\infty$. From this, we can take $a(\log n) = (\log
n)^{(v+1)/(v+2)}$, $b(\log n)\equiv0$ and $k(x) = u x^{v+2} / \{
(v+1)(v+2)\}$, thus the limiting survivor function is $\exp[-u x^{v+2}/
\{(v+1)(v+2)\}]$. Reversing the conditioning would lead to the
analogous result with $s,t$ in place of $u,v$.

A particular example of $V$ is the logistic dependence structure
(Tawn \cite{Tawn88}), which has $V(x,y) = (x^{-1/\alpha} +
y^{-1/\alpha})^\alpha$,
for $\alpha\in(0,1]$, and thus $\kappa(\beta,\gamma) = (\beta^{1/\alpha
} + \gamma^{1/\alpha})^\alpha$. Either directly or via the preceding
calculations we find for this example that $b(\log n) \equiv0$,
$a(\log n) = (\log n)^{1-\alpha}$ and $k(x) = \alpha x^{1/\alpha}$. In
this example, $\kappa$ is given by the $L_{1/\alpha}$ norm and provides
any $L$-norm between the $L_1$ and $L_\infty$ bounds given in
Property \ref{property4}.
\end{example}
%
\begin{example}[(Morgenstern)]
\label{secMorgenstern}
\[
C(u,v) = uv\bigl\{1+\alpha(1-u) (1-v)\bigr\};\qquad \alpha\in[-1,1].
\]
The expansion in exponential margins is
\begin{eqnarray*}
\Prob(X_E>\beta\log n,Y_E>\gamma\log n) &=& (1+
\alpha)n^{-(\beta
+\gamma)} - \alpha n^{-(2\beta+\gamma)}\\
&&{}- \alpha n^{-(\beta+2\gamma)}
+\alpha
n^{-(2\beta+2\gamma)}.
\end{eqnarray*}
This example highlights easily why the first order theory upon which we
focus is insufficient to make full connections to the conditional limit
theory described in Section \ref{secConditionedLimitTheory}. Here,
$L(n;\beta,\gamma) = (1+\alpha) - \alpha n^{-\beta}- \alpha
n^{-\gamma}
+\alpha n^{-(\beta+\gamma)}$, which collapses to 1 when $\beta$ or
$\gamma\rightarrow0$, but is in general different from 1 in its
leading order term.
\end{example}
%
\begin{example}[(Bivariate extreme value distribution)]
\label{secBivariateExtremeValueDistribution}
\[
C(u,v) = \exp\bigl[-V\{-1/\log u,-1/\log v\}\bigr],
\]
where $V$ is the exponent function. The joint survivor function in
exponential margins is
\begin{eqnarray*}
&&
\Prob(X_E>\beta\log n,Y_E>\gamma\log n) \\
&&\quad=
n^{-\beta} + n^{-\gamma} -1 \\
&&\qquad{}+ \exp\bigl[-V\bigl\{-1/\log
\bigl(1-n^{-\beta}\bigr), -1/\log\bigl(1-n^{-\gamma}\bigr)\bigr\}
\bigr].
\end{eqnarray*}
\end{example}
%
\begin{example}[(Lower joint tail of the Clayton distribution)]
\label{secLowerJointTailOfTheClaytonDistribution}
\[
C(u,v) = u+v-1 + \bigl\{(1-u)^{-1/\alpha} + (1-v)^{-1/\alpha} - 1 \bigr
\}^{-\alpha};\qquad \alpha>0.
\]
For the expansion in exponential margins, without loss of generality
assume $\gamma= \max\{\beta,\gamma\}$, then as $n\rightarrow\infty$
\begin{eqnarray*}
\Prob(X_E>\beta\log n,Y_E>\gamma\log n) & = & \bigl
\{n^{\beta/\alpha} + n^{\gamma/\alpha} - 1\bigr\}^{-\alpha}
\\
&=& n^{-\gamma}\bigl\{1+n^{(\beta-\gamma)/\alpha}-n^{-\gamma/\alpha}\bigr
\}^{-\alpha} \\
&=& n^{-\gamma} + \RMo\bigl(n^{-\gamma}\bigr).
\end{eqnarray*}
\end{example}
%
\begin{example}[(Trivariate example)]
\label{secTrivariateExample}
Let $T,U,V,W$ be mutually independent following standard Pareto
distributions, and define $X=\max\{T,U\}$, $Y=\max\{U,V\}$, $Z=\max\{
V,W\}$. Then in the obvious notation we have $\chi(X,Y,Z) = 0$, $\chi
(X,Y)>0$, $\chi(Y,Z)>0$, $\chi(X,Z)=0$. The identical marginal
distributions of $X,Y,Z$ are $\Prob(X\leq x) =:G(x)= (1-1/x)^2 =
1-2/x+1/x^2$. Transforming to standard Pareto margins, $X_P = 1/\{
1-G(X)\}$, so that as $n\rightarrow\infty$, $\Prob(X_P>n^\beta) =
\Prob
\{X>G^{-1}(1-n^{-\beta})\}= \Prob\{X> 2 n^{\beta} +\RMO(1)\}$. Thus for
all sufficiently large $n$, we can bound the joint survivor probability,
\begin{eqnarray*}
\Prob\bigl(X>3n^\beta,Y>3n^\gamma,Z>3n^\delta\bigr)
&\leq&\Prob\bigl(X_P>n^\beta,Y_P>n^\gamma,Z_P>n^\delta
\bigr) \\
&\leq&\Prob\bigl(X>n^\beta,Y>n^\gamma,Z>n^\delta
\bigr)
\end{eqnarray*}
and hence $\kappa$ can be derived from the leading order term of
$\Prob
(X>kn^\beta, Y>kn^\gamma, Z>kn^\delta)$, $k>0$. The expression for this
probability may be found in Wadsworth \cite{Wadsworth12}, from which
we can derive
\[
\kappa(\beta,\gamma,\delta) = \cases{\displaystyle \gamma+ \min\{\beta,\gamma,\delta
\},&\quad $
\gamma=\max\{\beta,\gamma,\delta\}$,
\cr
\beta+\delta, &\quad otherwise. }
\]
Note that this reduces to the correct two-dimensional $\kappa$
functions when any of $\beta,\gamma,\delta$ are set to~0; that is,
$\max
\{\gamma,\delta\}$, $\beta+\delta$ and $\max\{\beta,\gamma\}$,
respectively.
\end{example}


\section{Statistical methodology}
\label{secStatisticalMethodology}
We propose one possible approach for performing statistical inference
on extreme set probabilities, motivated by the theory of preceding
sections. The method builds upon the connections between $\lambda$ and
Pickands' dependence function, and is outlined in Section \ref
{secInferentialApproach}. Further methodological details and a
comparison of new and existing methodology are presented in
Section \ref
{secEstimationAndDiagnostics}. In the implementation, we assume that
the data have been transformed into exponential margins. This can be
achieved from arbitrary marginal distributions by a transformation such
as that described in Coles and Tawn \cite{ColesTawn91}.


\subsection{Inferential approach}
\label{secInferentialApproach}
We assume a sample of $m$ bivariate random vectors, and consider
estimation of probabilities that $(X_E,Y_E)$ lies in sets of the form
$(\omega\log n,\infty)\times((1-\omega)\log n, \infty)$. By
representation (\ref{A1}), for each $\omega\in[0,1]$,
%
\begin{equation}\label{Inf1}
\Prob\biggl(\min\biggl\{\frac{X_E}{\omega},\frac{Y_E}{1-\omega
} \biggr\} >\log n
+ \log t \biggr) = L^*(nt;\omega) (nt)^{-\lambda(\omega)}
\end{equation}
as $n\rightarrow\infty$, along with the conditional probability for $t>1$,
%
\begin{eqnarray}\label{Inf2}
&&\Prob\biggl(\min\biggl\{\frac{X_E}{\omega},\frac{Y_E}{1-\omega
} \biggr\} >\log n
+ \log t \Big| \min\biggl\{\frac{X_E}{\omega},\frac
{Y_E}{1-\omega} \biggr\}>\log n
\biggr)\nonumber\\[-8pt]\\[-8pt]
&&\quad = \frac
{L^*(nt;\omega
)}{L^*(n;\omega)}t^{-\lambda(\omega)} \rightarrow
t^{-\lambda(\omega)}.\nonumber
\end{eqnarray}
Since the expression in (\ref{Inf2}) has a regularly varying tail with
positive index in Pareto margins, one can use the Hill estimator of the
tail index (Hill \cite{Hill75}) as an estimator for $\lambda(\omega
)$. To
tidy notation, write $u_n=\log n$ and $v=\log t$. Some consequences of
the assumption on the ratio of slowly varying functions are explored at
the end of the section. Asymptotic consistency of the Hill estimator
requires that the number of exceedances, $k_{m,n}$, of the threshold
$u_n$ satisfies $k_{m,n}\to\infty$ as $m\to\infty$, but that
$k_{m,n}/m \to0$. Conditions for asymptotic normality are studied in
Haeusler and Teugels \cite{HaeuslerTeugels85} and de Haan and Resnick
\cite{deHaanResnick98}, for example.

Similar estimation procedures have been exploited previously in
methodology for extreme value problems: for max-stable distributions,
Pickands \cite{Pickands81} proposed an analogous estimator for the dependence
function $A$ when dealing with componentwise maxima data, see also
Hall and Tajvidi \cite{HallTajvidi00}. For the asymptotically
independent case,
Ledford and Tawn \cite{LedfordTawn96,LedfordTawn97} used a similar
scheme on the diagonal
$\omega=1/2$ to estimate the coefficient of tail dependence $1/\eta=
2\lambda(1/2)$. Our methodology extends this in a natural way to any
ray $\omega\in[0,1]$ in exponential margins.


\subsection{Estimation and diagnostics}
\label{secEstimationAndDiagnostics}
\begin{figure}

\includegraphics{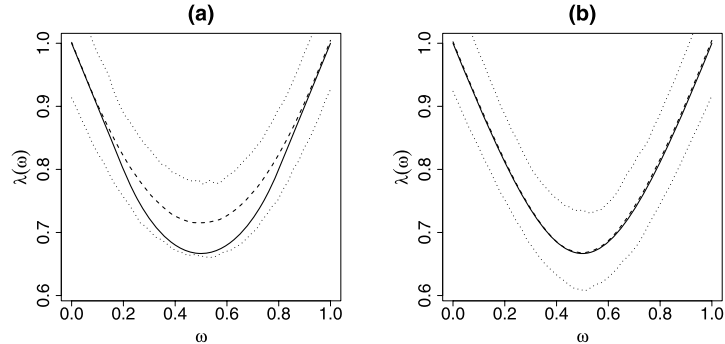}

\caption{True (solid line) and estimated mean (dashed line)
$\lambda(\omega)$, with 95\% confidence intervals (dotted lines) for
(a) the bivariate normal and (b) the inverted extreme value logistic
models. Mean and confidence intervals are based upon output from 500
simulated datasets of size 5000, with estimation based upon 10\% of the
sample.}
\label{figLambdaEst}
\end{figure}

We use the Hill estimator of the rate parameter, which has the simple
closed form of the reciprocal mean excess of the structure variable
$\min\{X_E/\omega, Y_E/(1-\omega)\}$ above the threshold $u_n$. Letting
$\hat{\lambda}(\omega)$ denote the value of the estimate, the final
estimate of the joint survivor probability is constructed as
\[
\hat{\Prob}\bigl\{X_E>\omega(u_n+v),
Y_E>(1-\omega) (u_n+v)\bigr\} = e^{-\hat
{\lambda}(\omega)v}
\tilde{\Prob}\bigl\{X_E>\omega u_n,
Y_E>(1-\omega) u_n\bigr\}
\]
with $\tilde{\Prob}$ representing the empirical exceedance probability.
In this simple methodology, that we employ in Section \ref
{secComparisonOfMethods}, there has been no attempt to impose any of
the constraints that $\lambda$ or $\kappa$ can satisfy, as identified
in Properties \ref{property1}--\ref{property4} and Proposition \ref
{prop3}. Imposition of some constraints could represent a
methodological extension if so desired. Figure \ref{figLambdaEst}(a)
and (b) display true and estimated $\lambda(\omega)$ functions under
this methodology for the examples considered in Section \ref
{secComparisonOfMethods}; the figures are discussed at the end of that
section.

To detect applicability of the theory, standard techniques such as
quantile--quantile plots can be used to assess the suitability of the
model. In addition, if the theory holds then for $m$ the sample size,
assumed large, we can write for $c>0$
\[
m\Prob\bigl\{X_E>c\omega\log m, Y_E>c(1-\omega) \log m
\bigr\} = L^*\bigl(m^c;\omega\bigr)m^{1-c\lambda(\omega)}.
\]
Since $L^*(m^c;\omega)$ is still slowly varying, it follows that on any
given ray $\omega$, plots of log number of points in the set $(c\omega
\log m,\infty)\times(c(1-\omega)\log m, \infty)$, for a variety of $c$,
should exhibit approximate linearity in $c$.


\subsection{Comparison of methods}
\label{secComparisonOfMethods}

The following simulation study compares the methodology of
Ledford and Tawn \cite{LedfordTawn97}, Heffernan and Tawn \cite
{HeffernanTawn04} and that of Section \ref
{secInferentialApproach} for estimating joint survivor probabilities
across a range of rays in exponential margins. The intention of the
comparison is to consider the properties of the different procedures
rather than to examine the evidence for an `optimal' procedure, since
we consider only simple implementations.

\begin{figure}

\includegraphics{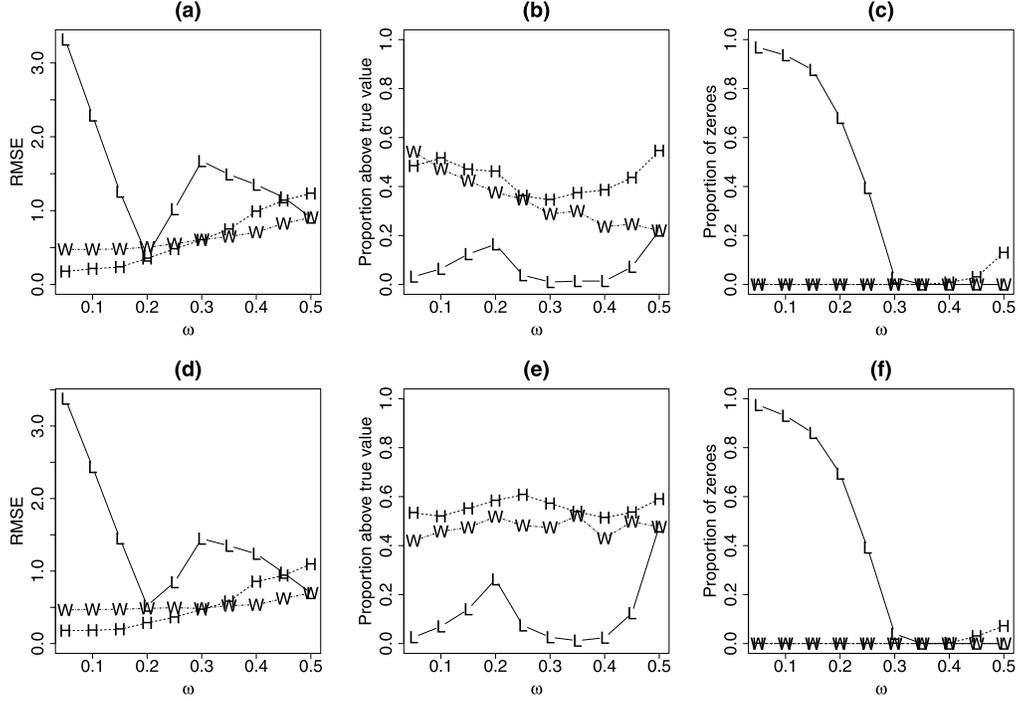}

\caption{(a): RMSE of non-zero log probabilities, (b):
proportion of estimates exceeding the true probability, and (c):
proportion of zero probability estimates versus $\omega$ for the
bivariate normal distribution with $\rho=0.5$; (d), (e), (f): as (a),
(b), (c) but for the inverted bivariate extreme value distribution with
logistic dependence structure and dependence parameter $\alpha=-\log
0.75/\log2$. Plot characters `H', `L', `W' refer to the methodology of
Heffernan and Tawn \cite{HeffernanTawn04}, Ledford and Tawn \cite
{LedfordTawn97} and the current authors,
respectively.}
\label{figSimResults}
\end{figure}

We examine two types of data which show the most flexible range of tail
behaviour amongst identified asymptotically independent distributions:
the bivariate normal distribution, and the inverted bivariate extreme
value distribution with logistic dependence structure, see Section \ref
{secConnectionsToExistingTheory}. For the results displayed below, the
parameters of the distribution were set such that $\eta=0.75$ in both
cases. For each type of data, we generated 5000 pairs and used $10\%$
of the data for estimation purposes. This process was repeated 500
times. We wish to assess the performance of the method across the whole
of the positive quadrant and thus we consider joint survivor sets
defined by pseudo-angles $\omega=0.5,0.45,\ldots,0.05$, fixing the $y$
coordinate at $1.5\log5000$ and multiplying by $\omega/(1-\omega)$ to
attain the $x$ coordinate. Our results are displayed in Figure \ref
{figSimResults} as root mean squared error (RMSE) of the non-zero log
probabilities; proportion of estimates exceeding the true value; and
proportion of probabilities estimated as exactly zero, all versus
$\omega$.

Overall, the methodology, similarly to the theory, provides something
of a compromise between that of Ledford and Tawn \cite{LedfordTawn97}
and of
Heffernan and Tawn \cite{HeffernanTawn04}. One advantage of the new
results is that one should
rarely be in a position where a probability is estimated as identically
zero, which we have argued can be a drawback to the existing
nonparametric methodology.

For $\omega= 0.5$, using the Hill estimator for $\eta$, the results
from the methodology of Ledford and Tawn \cite{LedfordTawn97} and the
current paper are
the same, since we have $1/\eta= \kappa(1,1) = 2\lambda(1/2)$. Away
from the diagonal however, the results of the Ledford and Tawn \cite
{LedfordTawn97}
model, which rely on estimating empirical probabilities through (\ref
{LText}) with $v+\log n + A=(1.5\log5000,\infty)\times(\{\omega
/(1-\omega)\}1.5\log5000,\infty)$, deteriorate greatly. This reflects
the fact that the theory underlying the methodology is designed for
joint tails, where both variables grow at the same rate. For angles
$\omega<0.3$ in this example, we observe a rapid rise in the number of
zero probability estimates, see Figure \ref{figSimResults}(c)
and (f). In any real application, this could be a serious problem as
only a single sample is available. The problem of having an empirical
probability of zero in a joint survivor set could in theory be
partially solved by noting that for any angle $\omega_1<\omega_2$ the
latter set is nested in the former, and one could sum the zero and
non-zero component. This however provides rather poor estimates, and so
we do not adopt this strategy for the reported results.

The methodology of Heffernan and Tawn \cite{HeffernanTawn04} shows
steady improvements in
RMSE as the angle $\omega$ decreases, and the joint survivor set moves
closer to being a marginal probability. This is owing to the form of
the limit theory involved. Recall equations (\ref{cevht}) and (\ref
{cevhr}). Let $\hat{a}, \hat{b}$ denote the estimates of the
normalizing functions, and define the $i$th component of a residual
vector by $Z_i:=\{X_{E,i}-\hat{b}(Y_{E,i})\}/\hat{a}(Y_{E,i})$. Making
$j=1,\ldots,r$ new draws $Y_{E,j}^*$ from the exponential distribution,
conditionally upon being larger than $(1-\omega) u_n$, and $r$ draws
from the empirical distribution function of the residuals, the required
probability can be estimated by
%
\begin{equation}\label{htest}
\Prob\bigl\{Y_{E}>(1-\omega) u_n\bigr\}
\frac{1}{r}\sum_{j=1}^r I\bigl\{
\hat{a}\bigl(Y_{E,j}^*\bigr)Z_j+\hat{b}
\bigl(Y_{E,j}^*\bigr)>\omega u_n |Y_{E,j}^*>(1-
\omega) u_n\bigr\}.
\end{equation}
When $(1-\omega) u_n$ is large but $\omega u_n$ small, it is the
marginal probability $\Prob\{Y_{E}>(1-\omega) u_n\}$ which is important
for accurate estimation. This ensures increasing precision of
estimation near the axes. The reliance of the estimator (\ref{htest})
upon an empirical residual vector highlights how zero probability
estimates are a potential problem, since only certain trajectories of
the new variable $\hat{a}(Y_{E,j}^*)Z_j+\hat{b}(Y_{E,j}^*)$ are
possible. This problem occurs only rarely here, though some evidence of
it for angles $\omega$ near 0.5 is apparent from Figure
\ref{figSimResults}(c) and (f), but the issue intensifies as the sample
size decreases or dimension increases. For an additional discussion,
see Peng and Qi \cite{PengQi04}.

The simple estimator we propose provides consistently reasonable
estimates across the whole range of angles $\omega$. In further
simulation results, not reported, we found that the RMSE can be
improved close to the axes ($\omega$ near 0 or 1), to a level similar
to the Heffernan and Tawn \cite{HeffernanTawn04} methodology, by
replacing the random
exponential sample with rank-transformed data; that is, the $i$th
largest value $X_{E,(i)}$ is replaced by $-\log\{i/(n+1)\}$.

Figure \ref{figLambdaEst}, along with Figure \ref{figSimResults}(b)
and (e) illustrate the effect that the non-constant slowly varying
function of the normal distribution has in comparison to the unit
slowly varying function of the inverted bivariate extreme value
distribution. Figure \ref{figSimResults}(b) shows a small bias on the
diagonal $\omega=0.5$ under the proposed estimation method, which
gradually decreases as the angle moves towards the axes. The true
conditional excess probability is
\[
\frac{L^*(nt;\omega)}{L^*(n;\omega)} t^{-\lambda(\omega)} = t^{-[\lambda
(\omega)-\{\log L^*(nt;\omega)-\log L^*(n;\omega)\}/\log t]}.
\]
For the normal distribution, estimation of $\lambda$ will therefore
exhibit some bias which is more pronounced in the range of $\omega$
such that $L_1(n;\omega,1-\omega)$, as identified in equation (\ref
{L1normal}), is non-constant; this is observed in Figure
\ref{figLambdaEst}(a). This bias can be seen to decay at rate $\RMO(1/\log
n)$, and thus will persist for nearly all practical sample sizes.
However, this fact ensures that the theoretical bias in $\lambda$
estimation is not as severe from a practical perspective, since it
remains largely unchanged over the range for which we are likely to be
interested in extrapolation. Figures \ref{figLambdaEst}(b) and
\ref{figSimResults}(e) show that the estimation of $\lambda$ for the
inverted bivariate extreme value distribution is unbiased for all
angles, as should be expected given that $L\equiv1$.


\section{Discussion}
\label{secDiscussion}
The theory presented has provided a characterization of tail
probability decay rates for a wide variety of multivariate
distributions. The results apply both in the case of asymptotically
independent and dependent random vectors. However, richer theoretical
representations are gained in the case of asymptotic independence in
two dimensions, or an analogous notion of weak joint tail dependence in
three or more dimensions, as defined in Section \ref
{secAsymptoticDependenceAndIndependence}. The theory developed
strongly mirrors the existing representation for asymptotically
dependent distributions, however it is exploited in the
characterization of tail probability decay rates rather than the
characterization of limit distributions. In particular, we find that
the homogeneous function $\kappa$ and angular dependence function
$\lambda$ play roles under weak joint tail dependence which are very
similar to the roles of the exponent function $V$ and Pickands'
dependence function $A$ of the classical theory for asymptotic
dependence. Multivariate extreme value theory is most often approached
from a perspective of multivariate regular variation. We have
demonstrated how generalizing multivariate regular variation of
$(X_P,Y_P)$ to that of $(X_P^{1/\omega},Y_P^{1/(1-\omega)})$ provides a
representation which encodes dependence information under both
asymptotic dependence and asymptotic independence.

From a statistical perspective, we have demonstrated that a simple
methodological procedure can provide very reasonable estimates of
bivariate joint survivor probabilities across the whole of the positive
quadrant. In addition, we have argued that in providing a theoretical
and methodological compromise between the approaches of
Ledford and Tawn \cite{LedfordTawn97} and Heffernan and Tawn~\cite
{HeffernanTawn04}, we can overcome some of
the methodological difficulties which can surface in either of these
two cases. In statistical applications, it is difficult to be certain
whether data are asymptotically independent or dependent, and to apply
the proposed methodology we need make no distinction.

The theory upon which we have focused is particularly suited to the
extrapolation of joint survivor sets, where only a single ray $\omega$
is implicated. Estimation of probabilities of more general sets can be
achieved in theory through exploitation of limit (\ref{Condconv2}).
However, rates of convergence to such limits may be slow. There is
clearly also a choice of limit distributions in this case, and
investigations into the most suitable trajectory of extrapolation for a
given set of interest, in terms of optimizing the rate of convergence
to the limit, remains an open avenue for further research.

\section*{Acknowledgements}

This research was conducted whilst JLW was based at Lancaster
University, and funding from the EPSRC and Shell Research through a
CASE studentship is gratefully acknowledged. We are grateful to Anthony
Ledford for discussions surrounding the bias that occurs in the
bivariate normal tail decay rate estimation, and to the referees and
Associate Editor for very helpful comments and suggestions which have
improved the paper substantially.


%

\printhistory

\end{document}